\documentclass[11pt]{amsart}
\usepackage{geometry}                
\usepackage{graphicx}
\usepackage{amssymb}
\usepackage{epstopdf}
\usepackage[all]{xy}
\usepackage{psfrag}
\usepackage[T1]{fontenc}
\usepackage[latin1]{inputenc}
\usepackage{a4}
\usepackage{graphics}
\usepackage{amsthm}
\usepackage{amssymb,graphicx,psfrag,lscape}
\usepackage{amssymb,epsf,colordvi,verbatim,graphics}
\usepackage{graphics,graphicx,psfrag,lscape}
\usepackage{epsf,epsfig,amsmath}
\usepackage{amssymb,latexsym}
\usepackage{lscape}

\newtheorem{thm}{Theorem}[section]

\newtheorem{pro}[thm]{Proposition}
\newtheorem{lem}[thm]{Lemma}

\newtheorem{cor}[thm]{Corollary}

\theoremstyle{definition}

\theoremstyle{remark}
\newtheorem*{rem}{Remark}

\title{A note on  Fontaine theory using different Lubin-Tate groups}
\medskip

\author{ Bruno Chiarellotto, Francesco Esposito}

\begin{document}
\maketitle
{\it Abstract.}
The starting point of Fontaine theory is the possibility of translating the study of   a $p$-adic representation of the absolute Galois group of a finite extension $K$ of $\mathbb{Q}_{p}$ into the investigation of a $(\varphi,\Gamma)$-module. This is done by decomposing the Galois group along  a totally ramified extension of $K$, via the theory of the field of norms:  the extension used is obtained by means of the cyclotomic tower which, in turn, is associated to the multiplicative Lubin-Tate group.
It is known that one can insert different   Lubin-Tate groups  into the  ''Fontaine theory" machine to obtain  equivalences with new categories of   $(\varphi,\Gamma)$-modules (here $\varphi$ may be iterated).  This article uses  only $(\varphi,\Gamma)$-theoretical terms to compare the different  $(\varphi,\Gamma)$ modules arising from various Lubin-Tate groups.

\let\thefootnote\relax

\section{Introduction and notations}
Let $K$ be a finite extension of $\mathbb{Q}_{p}$, $\overline{{K}}$ an algebraic closure of $K$, $G_{K}=\textrm{Gal}(\overline{K}/K)$ its absolute Galois group. In \cite{Fo}, Fontaine shows how the study of representations of $G_{K}$ over finite dimensional $\mathbb{Q}_{p}$-vector spaces is reduced to the study of an algebraic object, namely, the \'etale $(\varphi, \Gamma)-$module attached to it. 
In particular, he shows that there is a functor $D$ between the category ${\bf{Rep}}_{G_{K}}(\mathbb{Q}_{p})$ of $p$-adic representations of $G_{K}$ and the category ${\bf{\Phi}}^{et}_{\varphi,\Gamma}(\mathbb{B}_{K})$ of \'etale $(\varphi,\Gamma)$-modules over the discretely valued field $\mathbb{B}_{K}$, giving an equivalence of categories.
 The functor $D$ is defined as follows: $D(V)=(V\otimes_{\mathbb{Q}_{p}}\mathbb{B})^{H_{K}}$, where $\mathbb{B}$ is the completion of an unramified closure of $\mathbb{B}_{K}$; the group $H_{K}$ is the kernel of the cyclotomic character $\chi_{p}:G_{K}\rightarrow\mathbb{Z}_{p}^{\times}$. 

The field $\mathbb{B}_{K}$ has as residue field the field of norms  $\mathbb{E}_{K}$ of the extension $K(\bigcup_{n}\zeta_{ p^{n}})$ of $K$, obtained from $K$ by adding the $p^{n}$-th roots of 1 for every $n$.  The field $\mathbb{B}_{K}$ is endowed with the action of $\Gamma_{K}=\textrm{Im}(\chi_{p})$ and a commuting action of the Frobenius which lift those on $\mathbb{E}_{K}$. 

The $p^{n}$-th roots of unity are the $p^{n}$-th torsion points of the multiplicative Lubin-Tate formal group on $\mathbb{Q}_{p}$, associated to the uniformizer $p$. The cyclotomic tower is exactly the tower associated by Lubin-Tate theory to the multiplicative Lubin-Tate formal group. It is natural to try to carry out the theory for an arbitrary Lubin-Tate formal group $\mathcal{F}$ over a subfield $F$ of $K$, not only for ${\mathbb Q}_p$ (Partial results in this direction have been established by Kisin and Ren   \cite{Ki-R} and Fourquaux \cite{Fou}).

We now give   a more detailed description of our work along with that of \cite{Ki-R}.     Let $F$ be a subfield of $K$ ($K$ is a finite extension of ${\mathbb Q}_p$): $\mathcal{O}_{F}$ (resp. $\mathcal{O}_{K}$)  the ring of integers, $k_{F}$ (resp. $k_{K}$) the  residue field, $p^r$  (resp. $p^{rs})$the cardinality of $k_{F}$ (resp. $k_{K}$), and $\pi$  (resp. $\varpi$) a uniformizer of $F$ (resp. $K$).  We will denote by $F_0$ (resp. $K_0$) the maximal unramified extension of ${\mathbb Q}_p$ in $F$ (resp. in $K$). 

Let $F_{\pi}$ denote the maximal abelian totally ramified extension of $F$ associated to $\pi$ by Lubin-Tate theory; denote by $\chi_{\pi}: G_F \rightarrow {\mathcal O}_F^{\times}$ its associated character.  Let $X_{K}(K\cdot F_{\pi})= {\mathbb E}_{K,\pi}$ be the field 
of norms of the extension of $K$ which is the compositum of $K$ and $F_{\pi}$, and let its Galois group be denoted by 
$H_{K,\pi}$, so that $H_{K,\pi} = G_K \cap Ker( \chi_{\pi})$. The field $X_{K}(K\cdot F_{\pi})= {\mathbb E}_{K,\pi}$ is endowed with an action
 of $\Gamma_{K,\pi}=\textrm{Gal}(K\cdot F_{\pi}/K)$  and a commuting action of the Frobenius morphism $\varphi^{r}$ (relative to $k_{F}$). 
 Kisin and  Ren construct  (under the hypothesis of the  inclusion of $K$ in the compositum $K_0  \cdot  F_{\pi}$ \cite{Ki-R}2.1)    a  complete discrete valuation ring ,  $\mathcal{O}_{\mathcal{E}}$   (in our  notation $\mathbb{A}_{F,\pi}(K)$),  with residue field $X_{K}(K\cdot F_{\pi})={\mathbb E}_{K,\pi}$, unramified over $\mathcal{O}_{F}$. They are able to lift   the commuting actions of $\Gamma_{K,\pi}$ and $\varphi^{r}$ on $X_{K}(K\cdot F_{\pi})$ to commuting actions of $\Gamma_{K,\pi}$ and $\varphi^{r}$ on $\mathcal{O}_{\mathcal{E}} $  (   $\mathbb{A}_{K,\pi}(K)$ in our notation).
  Let ${\bf{Rep}}_{G_{K}}({\mathcal{O}}_{ F})$ be the category of representations of $G_{K}$ on finite free $\mathcal{O}_{F}$-modules and ${\bf{\Phi}}^{et}_{\varphi^{r},\Gamma_{K,\pi}}(\mathbb{A}_{K,\pi}(K))$  be the category of finite free $\mathcal{O}_{\mathcal{E}}(=\mathbb{A}_{K,\pi}(K))$-modules endowed with an \'etale action of $\varphi^{r}$ and a commuting action of $\Gamma_{K,\pi}$. 
  Kisin and Ren prove that the functor $V\mapsto(V\otimes\mathcal{O}_{\widehat{\mathcal{E}}^{nr}})^{H_{K,\pi}}$ is an equivalence of categories between ${\bf{Rep}}_{G_{K}}(\mathcal{O}_{F})$ and ${\bf{\Phi}}^{et}_{\varphi^{r},\Gamma_{K,\pi}}(\mathcal{O}_{\mathcal{E}})$, with as quasi-inverse the functor $M\mapsto (M\otimes\mathcal{O}_{\widehat{\mathcal{E}}^{nr}})^{\varphi^{r}=1}$, where $H_{K,\pi}$ is the absolute Galois group of $K\cdot F_{\pi}$.  In this  paper we want to give a more general formulation to these results:  we won't impose any condition on the extension $K/F$ and  we will treat  ${\mathcal O}_K$-representations. 
  Such results are perhaps known  in the math  common knowledge   (see, for example,  Fourquaux in  1.4.1  of  \cite{Fou}) , but we are not aware of any such explicit statement in the literature.  Namely, we prove that starting with a Lubin-Tate group $\mathcal{F}$ over $F$ associated to the uniformizer $\pi$ and for any finite extension $K$,  one may construct a complete discrete valuation ring $\mathbb{A}_{K,\pi}$ with residue field the field of norms $X_{K}(K\cdot F_{\pi})=\mathbb{E}_{K,\pi}$ and unramified over $\mathcal{O}_{K}$; moreover,  $\mathbb{A}_{K,\pi}$ is endowed with commuting actions of $G_{K}$ and $\varphi^{rs}$,  lifting those on $X_{K}(K\cdot F_{\pi})=\mathbb{E}_{K,\pi}$.
 We define as well, two functors $D_{\pi}$ and $V_{\pi}$. 
 We prove that the functor $D_{\pi}$ is an equivalence between ${\bf{Rep}}_{G_{K}}(\mathcal{O}_{K})$ and ${\bf{\Phi}}^{et}_{\varphi^{rs},\Gamma_{K,\pi}}(\mathbb{A}_{K,\pi})$; 
 and $V_{\pi}$ its quasi-inverse going from ${\bf{\Phi}}^{et}_{\varphi^{rs},\Gamma_{K,\pi}}(\mathbb{A}_{K,\pi})$ to ${\bf{Rep}}_{G_{K}}(\mathcal{O}_{K})$.
 
  Naturally one could have started with $F=K$ in the previous setting.  Then one would study  the Lubin-Tate group in $K$ associated to a choice $\varpi$ of a uniformizer of $K$.  Here the associated  totally ramified tower will define a factorization of $G_K$ by $H_{K,\varpi}$ and $\Gamma_{K,\varpi}$ 
  and the theory of the field of norms will yield  a   field ${\mathbb E}_{K,\varpi}$ together with a lifting ${\mathbb A}_{K, \varpi}$. Finally one will have a new equivalence between  ${\bf{Rep}}_{G_{K}}(\mathcal{O}_{K})$ and  
  ${\bf{\Phi}}^{et}_{\varphi^{rs},\Gamma_{K,\varpi}}(\mathbb{A}_{K,\varpi})$. 
  In particular these two categories of  $(\varphi,\Gamma)$-modules, namely ${\bf{\Phi}}^{et}_{\varphi^{rs},\Gamma_{K,\pi}}(\mathbb{A}_{K,\pi})$ and ${\bf{\Phi}}^{et}_{\varphi^{rs},\Gamma_{K,\varpi}}(\mathbb{A}_{K,\varpi})$ are equivalent: they are equivalent to the same category of  representations.  In the second part of our present work we will give  a proof  of the   equivalence in purely $(\varphi,\Gamma)$-module theoretic terms. 
 
 Let us now review in more detail the content of the various sections. 
 In section 2, we deal with  generalities on  Lubin-Tate theory.   Namely we study  maximal abelian ramified extensions of $F$ and $K$ given by Lubin-Tate theories  and their relative position. In particular we characterize when the two extensions are nested one inside the other.
 In section 3, we study the fields of norms attached to  two maximal abelian ramified extensions of $F$ and $K$ via    Lubin-Tate groups of different heights.  We consider  their relative position inside $\widetilde{\mathbb{E}}= \lim {\mathbb C}_p$. Here, the main result is that the  two extensions overlap only on the field of constants.
 In section 4,  we study the equivalences in positive characteristic between  ${\bf{Rep}}_{G_{K}}(\mathbb{F}_{p^{rs}})$,
${\bf{\Phi}}^{et}_{\varphi,\Gamma_{K,\pi}}({\mathbb{E}_{K,\pi}})$ and  ${\bf{\Phi}}^{et}_{\varphi,\Gamma_{K,\varpi}}({\mathbb{E}_{K,\varpi}})$. In particular, we describe the equivalence between the last two of these categories without appealing to the category ${\bf{Rep}}_{G_{K}}(\mathbb{F}_{p^{rs}})$.
In section 5, we begin to lift the theory to characteristic 0. In particular, we  recall and adapt  to our present case the construction of Kisin and Ren of the equivalence between the category of representations of $G_{F}$ on finite free $\mathcal{O}_{F}-$modules, and the category of \'etale $(\varphi^{r},\Gamma_{F})$-modules on finite free  $ \mathbb{A}_{F,\pi}$ modules (note that $\Gamma_{F}$ is linked to the totally ramified tower over $F$ given by $\pi$)
In section 6, we construct a complete discrete valuation ring  $\mathbb{A}_{K,\pi}$ endowed with actions of $G_{K}$ and $\varphi^{rs}$ and prove that the category ${\bf{Rep}}_{G_{K}}(\mathcal{O}_{K})$ is equivalent to the category $\Phi^{et}_{\varphi^{rs}, \Gamma_{K,\pi}}(\mathbb{A}_{K,\pi})$ of \'etale $(\varphi^{rs},\Gamma_{K,\pi})$-modules on finite free $\mathbb{A}_{K,\pi}$-modules.
In section 7, we give our  main comparison theorem between two  categories of $(\varphi,\Gamma)$-modules. We express these equivalence  functors  with the help of a system of rings $\mathbb{A}_{\pi, \varpi,\bullet}$ endowed with partial Frobenius maps and an action of $G_{K}$.

{\bf Acknowledgment:} The research of the two authors has been supported by the Eccellenza Grant Cariparo "Differential Methods in Algebra, Geometry and Analysis". We thank Lionel Fourquaux  for useful conversations. We thank Frank Sullivan.

\section{Lubin-Tate theory of $F$ and $K$}
Lubin-Tate theory gives an explicit construction of the maximal abelian extension of a local field and of the norm residue map of local classfield theory (\cite{Ha}, \cite{Se}).

\begin{pro}
Let $F\subset K$ be an extension of local fields, with $K$ finite over $\mathbb{Q}_{p}$. Let $\pi \in \mathcal{O}_{F}$ be a uniformizer of $F$, $\varpi \in \mathcal{O}_{K}$ be a uniformizer of $K$, let $F_{\pi}$ and $K_{\varpi}$ be the Lubin-Tate extensions of $F$ and $K$ respectively, attached to $\pi$ and $\varpi$. Then

\begin{itemize}
\item[(i) ] $F_{\pi}\subset K_{\varpi}$ if and only if $N_{K/F}(\varpi)=\pi ^{s}$, where $s$ is the index of inertia of the extension $K/F$.

\item[(ii)] If $K/F$ is Galois, then $K_{\varpi}$ is Galois over $F$ if and only if $K/F$ is unramified and $\varpi \in F$.

\end{itemize}
\end{pro}

{\bf Proof.}

(i) By \cite{Se} or \cite{LT}, it follows that $\textrm{Gal}(F^{ab}/K_{\varpi}\cap F^{ab})$ is the completion of $N_{K/F}(<\varpi>)$, hence by Galois correspondence $F_{\pi}\subset K_{\varpi}$ if and only if $N_{K/F}(<\varpi>)\subset <\pi>$, which is easily seen to be equivalent to  $N_{K/F}(\varpi)=\pi ^{s}$.

(ii) Since $K/F$ is Galois, we have that $K^{ab}/F$ is Galois. By functoriality of the norm residue symbol, it follows that, if $\sigma \in \textrm{Gal}(\overline{F}/F)$, then $\textrm{Gal}(K^{ab}/K_{\varpi}^{\sigma})$ is the completion of $<\varpi ^{\sigma}>$.

~$\Box$\\ 

\begin{rem}
If $K/F$ is ramified, the norm map $N_{K/F}$ from $\mathcal{O}_{K}^{\times}$ to $\mathcal{O}_{F}^{\times}$ is not surjective. This implies that for a
 general  choice of $\pi$ there need not be  
 a $\varpi$ such that $F_{\pi}\subset K_{\varpi}$; in fact, one can find such $K_{\varpi}$ if and only if $K\cdot F_{\pi}$ is totally ramified over $K$.  The fact that the image of the norm map is of finite index in $\mathcal{O}_{F}$ corresponds to the fact that upon taking a finite unramified extension $K^{'}$ of $K$ it is always possible to find such a $\varpi$; the smallest possible choice of $K^{'}$ being the unramified closure of $K$ inside $K\cdot F_{\pi}$.
\end{rem}

\section{the fields of norms}
Having chosen uniformizers $\pi \in \mathcal{O}_{F}$ and $\varpi \in \mathcal{O}_{K}$ such that $F_{\pi}\subset K_{\varpi}$, let $\mathcal{F}_{\pi}$ be a Lubin-Tate group over $\mathcal{O}_{F}$ attached to the uniformizer $\pi$ (see \cite{LT}) and let $[\pi]$ be the multiplication by $\pi$ in  $\mathcal{F}_{\pi}$. Let $\Lambda _{\pi,n}$ be the submodule of $\pi ^{n}$-torsion elements in the maximal ideal of $\mathcal{O}_{\overline{F}}$, \textit{i.e.},  $\Lambda _{\pi,n}=\{x\in \overline{F}| v_{F}(x)>0 \, \textrm{and} \, [\pi ^{n}](x)=0\}$,
and let  $K_{\pi, n}$ denote the extension of $K$ generated by $\Lambda _{\pi,n}$. We may define the field of norms, denoted 
$X_{K}(K\cdot F_{\pi})$ in \cite{W}, of the extension $K\cdot F_{\pi}/K$
as follows:
$$X_{K}(K\cdot F_{\pi})= \lim_{\leftarrow} K_{\pi, n},$$
where the projective limit is taken according to the norm maps. 

In the same way, $\mathcal{F}_{\varpi}$ denotes a Lubin-Tate group over $\mathcal{O}_{K}$ attached to the uniformizer $\varpi$, $\Lambda _{\varpi,n}=\{x\in \overline{K}| v_{K}(x)>0 \, \textrm{and} \, [\varpi ^{n}](x)=0\}$ and $K_{\varpi, n}=K[\Lambda _{\varpi,n}]$. We may define, according to \cite{W},  the field of norms $X_{K}(K_{\varpi})$, of the extension $K_{\varpi}/K$,
as follows:

$$ X_{K}(K_{\varpi})= \lim_{\leftarrow} K_{\varpi, n},$$
where the projective limit is  again taken according to the norm maps.

In what follows, we will denote $\mathbb{E}_{K,\pi}$ the field of norms of the extension $K\cdot F_{\pi}/K$, and $\mathbb{E}_{K,\varpi}$ the field of norms of the extension $K_{\varpi}/K$.

 We recall some results of \cite{W}, and give the main results concerning the relation between $\mathbb{E}_{K,\pi}$ and $\mathbb{E}_{K,\varpi}$.

\begin{thm}[see \cite{W} theorem 2.1.3]\label{W1}
$\mathbb{E}_{K,\pi}$ and $\mathbb{E}_{K,\varpi}$ have the structure of complete discretely valued fields of characteristic $p$ with residue field is isomorphic to $k_{K}$, the residue field of $K$; i.e. $\mathbb{E}_{K,\pi}\cong k_{K}((x))$ and $\mathbb{E}_{K,\varpi}\cong k_{K}((y))$. Moreover, $\mathbb{E}_{K,\pi}$ and $\mathbb{E}_{K,\varpi}$, are endowed with a natural action of $\Gamma_{K,\varpi}=\Gamma_{K}=\textrm{Gal}(K_{\varpi}/K)\cong \mathcal{O}_{K}^{\times}$, which on $\mathbb{E}_{K,\pi}$ factors through its image $\Gamma_{K,\pi}\subset \Gamma_{F} \cong \mathcal{O}_{F}^{\times}$ by the norm map $N_{K/F}$. 
\end{thm}

Let $H_{K,\pi}$, $H_{K,\varpi}$ be respectively $\textrm{Gal}(\overline{K}/K\cdot F_{\pi})$ and $\textrm{Gal}(\overline{K}/K_{\varpi})$.  We observe that $H_{K,\pi}/H_{K,\varpi}$ is isomorphic to the kernel of the norm map from $\mathcal{O}_{K}^{\times}$ to $\mathcal{O}_{F}^{\times}$. Let $\mathbb{E}_{\pi}$ be a separable algebraic closure of $\mathbb{E}_{K,\pi}$ and $G_{\mathbb{E}_{K,\pi}}$ be the absolute Galois group of $\mathbb{E}_{K,\pi}$; we will use the anologous notation for $\mathbb{E}_{K,\varpi}$.

\begin{thm}[see \cite{W} corollary 3.2.3, proposition 3.4.1]\label{W2}
The group $G_{\mathbb{E}_{K,\pi}}$ is canonically isomorphic to $H_{K,\pi}$ and $G_{\mathbb{E}_{K,\varpi}}$ is canonically isomorphic to $H_{K,\varpi}$; moreover the field of norms $\mathbb{E}_{K,\varpi}$ is isomorphic to the field of norms of the extension $\mathbb{E}_{\pi}^{H_{K,\varpi}}/\mathbb{E}_{K,\pi}$ and this isomorphism respects the action of $H_{K,\pi}$ ( which acts through the quotient $H_{K,\pi}/H_{K,\varpi}$).
\end{thm}

Let $\mathbb{C}_{p}$ denote the completion of $\overline{F}$. Let us consider the field $\tilde{\mathbb{E}}$ defined as the projective limit of copies of $\mathbb{C}_{p}$ by the map $x\mapsto x^{p}$:

$$\tilde{\mathbb{E}}=\lim_{\leftarrow}\mathbb{C}_{p}.$$

The basic properties of $\tilde{\mathbb{E}}$ are summarized in the following result (see \cite{FW}, \cite{W}).

\begin{thm}\label{W3}
The field $\tilde{\mathbb{E}}$ is an algebraically closed complete valued field of characteristic $p$ endowed with an action of $\textrm{Gal}( {\overline {\mathbb{Q}}}_{p}/\mathbb{Q}_{p})$. Moreover, the fields of norms $\mathbb{E}_{K,\pi}$ and $\mathbb{E}_{K,\varpi}$ have a natural Galois-equivariant immersion in $\tilde{\mathbb{E}}$ and the separable algebraic closures $\mathbb{E}_{\pi}$ and $\mathbb{E}_{\varpi}$ of $\mathbb{E}_{K,\pi}$ and $\mathbb{E}_{K,\varpi}$ respectively, are dense inside $\tilde{\mathbb{E}}$.
\end{thm}

A consequence of Theorem  \ref{W3} and the Theorem of Ax-Sen-Tate is the following.

\begin{cor}\label{W4}
The fixed fields $\tilde{\mathbb{E}}_{K,\pi}=\tilde{\mathbb{E}}^{H_{K,\pi}}$ and $\tilde{\mathbb{E}}_{K,\varpi}=\tilde{\mathbb{E}}^{H_{K,\varpi}}$ are respectively the completions of the perfections of $\mathbb{E}_{K,\pi}$ and $\mathbb{E}_{K,\varpi}$.
\end{cor}

\begin{rem}
Since $H_{K,\varpi}\subset H_{K,\pi}$ one has $\tilde{\mathbb{E}}_{K,\pi}\subset \tilde{\mathbb{E}}_{K,\varpi}$, \textit{i.e.}, the completions of the 
perfections 
 of the fields of norms are included one in the other in a Galois equivariant way. This does not continue to hold at the level of the fields of norms, even if $\mathbb{E}_{K,\pi}$ and $\mathbb{E}_{K,\varpi}$ are both abstractly isomorphic to the field $k_{K}((t))$ for an indeterminate $t$.
\end{rem}

We now show that, if one takes into account the action of $H_{K,\pi}$, the only possible common subfield of $\mathbb{E}_{K,\pi}$ and $\mathbb{E}_{K,\varpi}$ is their residue field $k_{K}$.

\begin{lem}\label{lem0}
If $F\neq K$, we have $ \mathbb{E}_{K,\varpi}^{H_{K,\pi}}=k_{K}$.
\end{lem}

{\bf Proof.}
As stated in Theorem \ref{W2}, the field $\mathbb{E}_{K,\varpi}$ may be identified with the field of norms of the extension $\mathbb{E}_{\pi}^{H_{K,\varpi}}/\mathbb{E}_{K,\pi}$, which is the projective limit along the inverse system of finite separable subextensions of $\mathbb{E}_{\pi}^{H_{K,\varpi}}/\mathbb{E}_{K,\pi}$ by the norm maps. The extension $\mathbb{E}_{\pi}^{H_{K,\varpi}}/\mathbb{E}_{K,\pi}$ is Galois with Galois group $H_{K,\pi}/H_{K,\varpi}$. This group is isomorphic to the kernel of the norm map $N_{K/F}:\mathcal{O}_{K}^{\times}\rightarrow \mathcal{O}_{F}^{\times}$, hence it is abelian. This implies that any finite separable subextension of  $\mathbb{E}_{\pi}^{H_{K,\varpi}}/\mathbb{E}_{K,\pi}$ is necessarily Galois; which implies that the group $H_{K,\pi}$ acts on the projective limit componentwise.  
Thus the elements fixed by 
$H_{K,\pi}$ have all components in $\mathbb{E}_{K,\pi}=\mathbb{E}_{\pi}^{H_{K,\pi}}$ and the norm maps coincide with the elevation to the degree of the extension on $\mathbb{E}_{K,\pi}$. 

Observe that the group $\mathcal{O}_{K}^{\times}$ is a $p$-adic Lie group of dimension $[K:\mathbb{Q}_{p}]$ and $\mathcal{O}_{F}^{\times}$ is a $p$-adic Lie group of dimension $[F:\mathbb{Q}_{p}]$. The norm map $N:\mathcal{O}_{K}^{\times}\rightarrow \mathcal{O}_{F}^{\times}$ is an open morphism of $p$-adic Lie groups. This implies that the kernel of $N$ is a $p$-adic Lie group of dimension $[K:F]$. Hence, except when $F=K$, the group $H_{K,\pi}/H_{K,\varpi}$ is infinite.

We use $U_{1}$ to denote the subgroup of $\mathcal{O}_{K}^{\times}$ of those elements congruent to 1 modulo the maximal ideal of $\mathcal{O}_{K}^{\times}$, and  $V_{1}$ to indicate the subgroup of $\mathcal{O}_{F}^{\times}$ of those elements congruent to 1 modulo the maximal ideal of $\mathcal{O}_{F}^{\times}$. It is well known that $U_{1}$ and $V_{1}$ are pro-$p$-groups of finite index in $\mathcal{O}_{K}^{\times}$ and $\mathcal{O}_{F}^{\times}$ respectively, and that the map $N$ sends $U_{1}$ to $V_{1}$. This implies that the kernel of $N$ contains a pro-$p$-subgroup of finite index, namely the intersection of the kernel of $N$ with $U_{1}$.

This proves that, in the case $F\neq K$, the group $H_{K,\pi}/H_{K,\varpi}$ has a finite index subgroup which is an infinite pro-$p$-group.

Since the inverse limit does not change upon restricting to a cofinal set of indices  (and in particular to  a finite extension of the base field ${\mathbb E}_{K,\pi}$)  : one may take the cofinal set of extensions which is made up of all finite extensions of the extension of $\mathbb{E}_{K,\pi}$ which corresponds to this finite index infinite pro-$p$-group. In this way the extensions are all $p$-extensions of one another.   So an element of the field of norms of the extension $\mathbb{E}_{\pi}^{H_{K,\varpi}}/\mathbb{E}_{K,\pi}$ is seen to be an infinite sequence $(x_{i})$ of elements of $\mathbb{E}_{K,\pi}$ coherent under the Frobenius morphism: $x_{i}=x_{i+1}^{p}$.
So, $\forall i,j \in \mathbb{N}, x_{i}\in (\mathbb{E}_{K,\pi})^{p^{j}}$. But the intersection of all $(\mathbb{E}_{K,\pi})^{p^{j}}$ is the field of constants $k_{K}$ of $\mathbb{E}_{K,\pi}$. Hence  the elements of the field of norms of the extension  $\mathbb{E}_{\pi}^{H_{K,\varpi}}/\mathbb{E}_{K,\pi}$ are exactly those given by taking coherent sequences of elements of $k_{K}$ under the Frobenius, and so the field of norms is isomorphic to $k_{K}$. ~$\Box$\\

\begin{cor}\label{cor2}
Suppose $F\neq K$.  Then inside $\tilde{\mathbb{E}}$ one has that  the following holds for the intersection of  $\mathbb{E}_{K,\varpi}$ and $\tilde{\mathbb{E}}_{K,\pi}$:

$$ \mathbb{E}_{K,\varpi} \cap \tilde{\mathbb{E}}_{K,\pi}=k_{K}.$$

The statement continues to hold even after replacing $\mathbb{E}_{K,\varpi}$ by its perfect closure inside $\tilde{\mathbb{E}}$. 
Moreover, if $x\in \tilde{\mathbb{E}}$ is a generating element of $\mathbb{E}_{K,\pi}$ over $k_{K}$ and $y\in\tilde{\mathbb{E}}$ is a generating element of $\mathbb{E}_{K,\varpi}$ over $k_{K}$, then $x$ and $y$ are algebraically independent  over $k_{K}$.
\end{cor}

{\bf Proof.} The first assertion is a consequence of Lemma \ref{lem0} because of the fact that $ \mathbb{E}_{K,\varpi} \cap \tilde{\mathbb{E}}_{K,\pi}$ is included in $\mathbb{E}_{K,\varpi}^{H_{K,\pi}}$. 
On the other hand, if $\alpha$ is an element in $\tilde{\mathbb{E}}$ which belongs to $\tilde{\mathbb{E}}_{K,\pi}$ and is purely inseparable over $\mathbb{E}_{K,\varpi}$, then a sufficiently high  $p$-power of $\alpha$ will be in $\mathbb{E}_{K,\varpi} \cap \tilde{\mathbb{E}}_{K,\pi}$. Hence $\alpha$ itself will be in $k_{K}$.

To prove the second, one may argue as follows.
If $x$ and $y$ were to satisfy an algebraic relation with coefficients in $k_{K}$, then $y$ would  be algebraic over $\tilde{\mathbb{E}}_{K,\pi}$ and the group of the continuous automorphisms of $\tilde{\mathbb{E}}$ which fix $y$ and $\tilde{\mathbb{E}}_{K,\pi}$ would have finite index in the group of the continuous automorphisms of $\tilde{\mathbb{E}}$ which fix $\tilde{\mathbb{E}}_{K,\pi}$. By Corollary \ref{W4}, the first group is $H_{K,\varpi}$ and the second is $H_{K,\pi}$. Since  $H_{K,\varpi}$ is not of finite index in $H_{K,\pi}$, $x$ and $y$ are algebraically independent over $k_{K}$.
~$\Box$\\

\begin{cor}
In the case $F\neq K$, there is no $G_{K}$-equivariant inclusion of $\mathbb{E}_{K,\pi}$ in $\mathbb{E}_{K,\varpi}$, nor is there one between their perfect closures.
\end{cor}

A basic consequence of Corollary \ref{cor2} is the following.

\begin{lem}\label{lem2}
In the case $F\neq K$, in $\tilde{\mathbb{E}}$, we have $\mathbb{E}_{\pi}\cap \mathbb{E}_{\varpi}=\overline{k}_{K}$, and the two extensions are algebraically independent over $\overline{k}_{K}$
\end{lem}

{\bf Proof.}
By Corollary \ref{cor2}, we know that $\mathbb{E}_{K,\pi}=k_{K}((x))$ and $\mathbb{E}_{K,\varpi}=k_{K}((y))$, with $x$ and $y$ algebraically independent. Suppose that $\alpha\in\mathbb{E}_{\pi}\cap\mathbb{E}_{\varpi}$. The element $\alpha$ has minimal polynomial $P(T)\in\mathbb{E}_{K,\pi}[T]$ over $\mathbb{E}_{K,\pi}$ and minimal polynomial $Q(T)\in\mathbb{E}_{K,\varpi}[T]$ over $\mathbb{E}_{K,\varpi}$. These polynomials will still be irreducible over the composite field of $\mathbb{E}_{K,\pi}$ and $\mathbb{E}_{K,\varpi}$ inside $\tilde{E}$ since $x$ and $y$ are algebraically independent variables over $k_{K}$. Hence the polynomials must coincide, and thus have coefficients in $k_{K}$.
$\Box$\\ 

\section{The system of fields $\mathbb{E}_{12,\bullet}$}

We first develop the theory in characteristic $p$. 
Let $W$ be a representation of $G_{K}$ over $\mathbb{F}_{p}$. To $W$ we may associate $D_{\pi}(W)=(W\otimes_{\mathbb{F}_{p}}\mathbb{E}_{\pi})^{H_{K,\pi}}$, 
which is a $(\varphi,\Gamma_{K,\pi})$-module over the field of norms
 $\mathbb{E}_{K,\pi}$. A fundamental result of \cite{Fo} is that the functor $D_{\pi}$ gives an equivalence of categories between 
 ${\bf{Rep}}_{G_{K}}(\mathbb{F}_{p})$ and 
$\bf{\Phi^{et}_{\varphi,\Gamma_{K,\pi}}}({\mathbb{E}_{K,\pi}})$
. Analogously, the functor $D_{\varpi}(W)=(W\otimes_{\mathbb{F}_{p}}\mathbb{E}_{\varpi})^{H_{K,\varpi}}$ gives an equivalence of categories between 
${\bf{Rep}}_{G_{K}}(\mathbb{F}_{p})$ and ${\bf{\Phi}}^{et}_{\varphi,\Gamma_{K,\varpi}}({\mathbb{E}_{K,\varpi}})$. 

We show that the argument in \cite{Fo} may be applied to prove that there is an equivalence between the category ${\bf{Rep}}_{G_{K}}(k_{K})$ and the category ${\bf{\Phi}}^{et}_{\varphi^{rs},\Gamma_{K,\pi}}({\mathbb{E}_{K,\pi}})$.

\begin{thm}\label{teor1}
The functor $D_{\pi}(W)=(W\otimes_{k_{K}}\mathbb{E}_{\pi})^{H_{K,\pi}}$ is an equivalence of categories between the category ${\bf{Rep}}_{G_{K}}(k_{K})$ and the category ${\bf{\Phi}}^{et}_{\varphi^{rs},\Gamma_{K,\pi}}({\mathbb{E}_{K,\pi}})$. The functor $V_{\pi}(D)=(D\otimes_{\mathbb{E}_{K,\pi}\mathbb{E}_{\pi}})^{\varphi^{rs}=1}$ is quasi-inverse to $D$.
\end{thm}
{\bf Proof.} Our argument is  an adaptation of Fontaine's proof  \cite{Fo} (see also \cite{Ki-R}).
There are natural maps 
$$D_{\pi}(W)\otimes_{\mathbb{E}_{K,\pi}}\mathbb{E}_{\pi}\rightarrow W\otimes_{k_{K}}\mathbb{E}_{\pi},$$
and
$$V_{\pi}(D)\otimes_{k_{K}}\mathbb{E}_{\pi}\rightarrow D\otimes_{\mathbb{E}_{K,\pi}}\mathbb{E}_{\pi}.$$
These natural maps commute with the action of $G_{K}$ and $\varphi^{rs}$, so to prove the statement of the theorem it is sufficient to prove that these two natural maps are isomorphisms. 

Since $H_{K,\pi}$ is the Galois group of the extension $\mathbb{E}_{\pi}$ over $\mathbb{E}_{K,\pi}$, the first map is an isomorphism by Hilbert's Theorem 90, i.e. $H^{1}(H_{K,\pi}, \textrm{GL}(n,\mathbb{E}_{\pi}))=\{ 1\}$. 

For the second map, observe that it is injective. In fact, a family of elements of $M_{\pi}(D)$ which is linearly independent over $k_{K}$ will remain linearly independent over $\mathbb{E}_{\pi}$. 
Indeed, let us take a family $\eta_{1},\ldots,\eta_{t}\in V_{\pi}(D)$ of elements which are independent over $k_{K}$ and such that any strict subfamily of this family remains linearly independent over $\mathbb{E}_{\pi}$. Suppose that the elements $x_{1},\ldots,x_{t}\in\mathbb{E}_{\pi}$ are all nonzero and such that $x_{1}\eta_{1}+\ldots+x_{t}\eta_{t}=0$. Dividing by $x_{1}$, we may suppose $x_{1}=1$. Now, applying $\varphi^{rs}$ to the relation we get $\eta_{1}+x_{2}^{p^{rs}}\eta_{2}+\ldots+x_{t}^{p^{rs}}\eta_{t}=0$; so , by subtracting one relation from the other one gets $(x_{2}^{p^{rs}}-x_{2})\eta_{2}+\ldots+(x_{t}^{p^{rs}}-x_{t})\eta_{t}=0$. This implies $x_{j}^{p^{rs}}-x_{j}=0$, i.e. all coefficients $x_{1},\ldots,x_{t}$ are in $k_{K}$.

Hence to prove the second map bijective, it is sufficient to prove that $V_{\pi}(D)$ has dimension over $k_{K}$ the dimension of $D$ over $\mathbb{E}_{K,\pi}$.   

 Let us choose a basis of the $\mathbb{E}_{K,\pi}$-vector space $D$, so that $D$ is identified with the $\mathbb{E}_{K,\pi}$-vector space $\mathbb{E}_{K,\pi}^{d}$. The action of $\varphi^{rs}$ is identified with an invertible matrix $(a_{j,l})$; i.e., 
$$
\varphi^{rs}x_{j}=\sum_{1\leq l\leq d} a_{j,l}x_{l}.
$$
The $k_{K}$-vector space $V_{\pi}(D)$ is then identified with the sub-$k_{K}$-vector space of $\mathbb{E}_{\pi}^{d}$ of the $d$-uples $(x_{j})_{1\leq j\leq d}$ verifying

$$
x^{p^{rs}}_{j}=\sum a_{j,l}x_{l}.
$$
So $V_{\pi}(D)$ is the set of $\mathbb{E}_{\pi}$-points of the \'etale $\mathbb{E}_{K,\pi}$-algebra $C=\mathbb{E}_{K,\pi}[X_{1}, \ldots , X_{d}]/(X^{p^{rs}}_{j}=\sum a_{j,l}X_{l})_{1\leq j\leq d}$.
Since $C$ has rank $p^{rsd}$, the dimension of $V_{\pi}(D)$ over $k_{K}$ is $d$. 
$\Box$\\ 

Naturally, the same arguments work in the case of $D_{\varpi}$ and $V_{\varpi}$. Thus we have also the following.
\begin{thm}\label{teor2}
The functor $D_{\varpi}(W)=(W\otimes_{k_{K}}\mathbb{E}_{\varpi})^{H_{K,\varpi}}$ is an equivalence of categories between the category ${\bf{Rep}}_{G_{K}}(k_{K})$ and the category ${\bf{\Phi}}^{et}_{\varphi^{rs},\Gamma_{K,\varpi}}({\mathbb{E}_{K,\varpi}})$. The functor $V_{\varpi}(D)=(D\otimes_{\mathbb{E}_{K,\varpi}\mathbb{E}_{\varpi}})^{\varphi^{rs}=1}$ is quasi-inverse to $D_{\varpi}$.
\end{thm}

{\bf Proof.}
As in \ref{teor1}.
$\Box$\\

We now wish to compare the modules $D_{\pi}(W)$ and $D_{\varpi}(W)$.

Let $\sigma$ be the $rs$-th power of the Frobenius automorphism of $\overline{k}_{K }$. We define $\mathbb{E}_{\pi,\varpi,i}$ as the field of fractions of the ring $\mathbb{E}_{\pi}\otimes_{\sigma^{i}}\mathbb{E}_{\varpi}$. By Lemma \ref{lem2}, the two fields are algebraically independent over $\overline{k}_{K}$, so the tensor product is a domain. The subscript $\sigma^{i}$ means that if $\lambda\in\overline{k}$, in  $\mathbb{E}_{\pi}\otimes_{\sigma^{i}}\mathbb{E}_{\varpi}$, we have $\lambda\otimes 1=1\otimes \sigma^{i}(\lambda)$. On every field $\mathbb{E}_{12,i}$ there is an action of $G_{K}$, acting diagonally on $\mathbb{E}_{\pi}$ and $\mathbb{E}_{\varpi}$. Moreover, the $rs$-th iterate $\varphi_{\pi}$ of the Frobenius morphism  of $\mathbb{E}_{\pi}$ induces maps $\varphi_{\pi}:\mathbb{E}_{\pi,\varpi,i}\rightarrow \mathbb{E}_{\pi, \varpi,i-1}$; the $rs$-th iterate $\varphi_{\varpi}$ of the Frobenius morphism  of $\mathbb{E}_{\varpi}$ induces maps $\varphi_{\varpi}:\mathbb{E}_{\pi,\varpi,i}\rightarrow \mathbb{E}_{\pi, \varpi,i+1}$. These maps commute with the action of $G_{K}$, and $\varphi_{\pi}\circ\varphi_{\varpi}=\varphi_{\varpi}\circ\varphi_{\pi}=\varphi^{rs}$ the $rs$-th iterate of the absolute Frobenius morphism. So we have a doubly projective system of fields with a $G_{K}$-action.

\begin{pro}\label{pro1}
The projective limit of the system $\mathbb{E}_{\pi, \varpi,\bullet}$ along $\varphi_{\pi}$ is $\mathbb{E}_{\varpi}$, and the morphisms $\varphi_{\varpi}$ are transformed in the $rs$-th iterate  $\varphi_{\varpi}$ of the Frobenius morphism of $\mathbb{E}_{\varpi}$. The symmetrical statement holds also.
\end{pro}
{\bf Proof.}
The field $\mathbb{E}_{\varpi}$ is the union of the finite separable extensions of $\mathbb{E}_{K,\varpi}$. These extensions are all of the type $k'((t))$ with $k'$ a finite extension of $k$ and $t$ a variable. So, any element in $\mathbb{E}_{\pi,\varpi,i}$ is inside the field of quotients of $\mathbb{E}_{\pi}\otimes_{k',\sigma^{i}}k'((t))$. This field is inside the field $\mathbb{E}_{\pi}((t))$; $\varphi_{\pi}$ acts only on ${\mathbb{E}}_{\pi}$, in the usual way. So, an element is in the projective limit of the system $\mathbb{E}_{\pi,\varpi ,i}$ if and only if it is a constant sequence of elements in $\mathbb{E}_{\varpi}$. Clearly this can be identified with $\mathbb{E}_{\varpi}$, and the morphism $\varphi_{\varpi}$ is identified with the usual $rs$-th iterate of the Frobenius morphism of $\mathbb{E}_{\varpi}$. The symmetrical statement is proved symmetrically.
$\Box$\\ 

We now define the two functors $\Phi$ and $\Psi$ between the categories  ${\bf{\Phi}}^{et}_{\varphi^{rs},\Gamma_{K,\pi}}({\mathbb{E}_{K,\pi}})$ and  ${\bf{\Phi}}^{et}_{\varphi^{rs},\Gamma_{K,\varpi}}({\mathbb{E}_{K,\varpi}})$ which give the equivalence of these two categories.

Let $D_{1}$ be in ${\bf{\Phi}}^{et}_{\varphi,\Gamma_{K,\pi}}({\mathbb{E}_{K,\pi}})$, we form the system $D_{1}\otimes_{\mathbb{E}_{K,\pi}}\mathbb{E}_{\pi, \varpi,\bullet}$. We make $\varphi_{\pi}$ and $G_{K}$ act diagonally on $D_{\pi}$ and $\mathbb{E}_{\pi, \varpi, \bullet}$. Analogously, if $D_{2}$ is in ${\bf{\Phi}}^{et}_{\varphi,\Gamma_{K,\varpi}}({\mathbb{E}_{K,\varpi}})$, then we may form the system $D_{2}\otimes_{\mathbb{E}_{K,\varpi}}\mathbb{E}_{\pi, \varpi,\bullet}$, making $\varphi_{\varpi}$ and $G_{K}$ act diagonally on $D_{2}$ and $\mathbb{E}_{\pi,\varpi,\bullet}$. Let $\Phi(D_{1})=(\textrm{projlim}_{\varphi_{\pi}}D_{1}\otimes_{\mathbb{E}_{K,\pi}}\mathbb{E}_{\pi, \varpi,\bullet})^{H_{K,\varpi}}$; it is endowed with an action of $\Gamma_{K,\varpi}$ and $\varphi_{\varpi}$. Analogously, let  $\Psi(D_{2})=(\textrm{projlim}_{\varphi_{\varpi}}D_{\varpi}\otimes_{\mathbb{E}_{K,\varpi}}\mathbb{E}_{\pi, \varpi,\bullet})^{H_{K,\pi}}$; it is endowed with an action of $\Gamma_{K,\pi}$ and $\varphi_{\pi}$.

\begin{pro}\label{pro2}
If $D_{1}$ is in ${\bf{\Phi}}^{et}_{\varphi^{rs},\Gamma_{K,\pi}}({\mathbb{E}_{K,\pi}})$, then $\Phi(D_{1})$ is in ${\bf{\Phi}}^{et}_{\varphi^{rs},\Gamma_{K,\varpi}}({\mathbb{E}_{K,\varpi}})$. 
If $D_{2}$ is in ${\bf{\Phi}}^{et}_{\varphi^{rs},\Gamma_{K,\varpi}}({\mathbb{E}_{K,\varpi}})$, then
 $\Psi(D_{2})$ is in ${\bf{\Phi}}^{et}_{\varphi^{rs},\Gamma_{K,\pi}}({\mathbb{E}_{K,\pi}})$. 
These functors are quasi-inverses of each other.
\end{pro}

{\bf Proof.}
We have $D_{1}\otimes_{\mathbb{E}_{K,\pi}}\mathbb{E}_{\pi,\varpi,i}=D_{1}\otimes_{\mathbb{E}_{K,\pi}}\mathbb{E}_{\pi}\otimes_{\mathbb{E}_{\pi}}\mathbb{E}_{\pi,\varpi,i}$. By Theorem \ref{teor2}, we have that $D_{1}$ comes from a representation $W$ of $G_{K}$ over $k_{K}$, and that $D_{1}\otimes_{\mathbb{E}_{K,\pi}}\mathbb{E}_{\pi}=W\otimes_{k_{K}}\mathbb{E}_{\pi}$.  So $D_{1}\otimes_{\mathbb{E}_{K,\pi}}\mathbb{E}_{\pi,\varpi,i}=W\otimes_{k_{K}}\mathbb{E}_{\pi,\varpi,i}$.
 Hence $\Phi(D_{1})=(\textrm{projlim}_{\varphi_{\pi}}D_{1}\otimes_{\mathbb{E}_{K,\pi}}\mathbb{E}_{\pi,\varpi,\bullet})^{H_{K,\varpi}}$ is equal to $(\lim_{\leftarrow,\varphi_{\pi}}W\otimes_{k_{K}}\mathbb{E}_{\pi, \varpi,\bullet})^{H_{K,\varpi}}$, which by Proposition \ref{pro1} is equal to $(W\otimes_{k_{K}}\mathbb{E}_{\varpi})^{H_{K,\varpi}}$ which is the $(\varphi^{rs},\Gamma_{K,\varpi})$-module \'etale associated to $W$. This shows that $\Phi(D_{\pi}(W))=D_{\varpi}(W)$. Analogously one proves that $\Psi(D_{\varpi}(W))=D_{\pi}(W)$. The two functors are thus quasi inverses and commute with the functors $D_{\pi}$ and $D_{\varpi}$.
$\Box$\\

\section{The ring $\mathbb{A}_{F,\pi}$}

In this paragraph we want to lift to characteristic zero the action of Galois and Frobenius on the field of norms.  In this paragraph nothing is new.  We  just adapt to our notations and setting  results presented in \cite{Ki-R} in the case $F=K$ ($L=K$ in notation of  \cite{Ki-R}). Let $F$ be a finite extension of $\mathbb{Q}_{p}$, $\mathcal{O}_{F}$ its ring of integers, $k_{F}$ its residue field, $p^{r}$ the cardinality of $k_{F}$, $\pi$ a uniformizer of $F$, $\overline{F}$ an algebraic closure of $F$ and $G_{F}$ its absolute Galois group. Lubin-Tate theory allows to construct a maximal abelian ramified extension $F_{\pi}$ of $F$ (\cite{LT}). Let $\chi_{\pi}:G_{F}\rightarrow \mathcal{O}_{F}^{\times}$ be the associated character. Let us consider the field of norms $X_{F}(F_{\pi})=\mathbb{E}_{F,\pi}$ of the extension $F_{\pi}$ of $F$ (see \cite{FW}, \cite{W} or \cite{Co}). It is a local field of characteristic $p$, with residue field $k_{F}$, such that its absolute Galois group $H_{F}=\textrm{Ker} \chi_{\pi}$ is naturally isomorphic to the absolute Galois group  of the field $F_{\pi}$. Moreover, $\mathbb{E}_{F,\pi}$ comes equipped with an action of the group $\Gamma_{F}=G_{F}/H_{F}\cong \mathcal{O}_{F}^{\times}$.

Let $\widetilde{\mathbb{E}}^{+}$ be the inverse limit of $\mathcal{O}_{\overline{F}}/p\mathcal{O}_{\overline{F}}$ with transition maps given by the Frobenius map. Its field of fractions, denoted by  $\widetilde{\mathbb{E}}$,  is a perfect field of characteristic $p$. Moreover there is a natural map $\mathbb{E}_{F,\pi}\rightarrow\widetilde{\mathbb{E}}$ (see \cite{W} 4.2) compatible with the Frobenius and the action of $G_{F}$ (through $\Gamma_{F}$). Let $\widetilde{\mathbb{A}}=W(\widetilde{\mathbb{E}})$ be the ring of Witt vectors with coefficients in the perfect field $\widetilde{\mathbb{E}}$ and $\widetilde{\mathbb{B}}=\textrm{Frac}(\widetilde{\mathbb{A}})$.
\begin{lem}
There exists a complete discrete valuation ring,  $\mathbb{A}_{F,\pi}$ with uniformizer $\pi$ and residue field $\mathbb{E}_{F,\pi}$, endowed with commuting actions of $\Gamma_{F}$ and $\varphi^{r}$, lifting the actions of $\Gamma_{F}$ and $\varphi^{r}$ on $\mathbb{E}_{F,\pi}$.
\end{lem}
{\bf Proof.}
This is the construction in \cite{Ki-R} (1.2 and 1.3). Let us recall briefly the main steps.  It is possible to find an element $ u \in \widetilde{\mathbb{A}} \otimes_{{\mathcal O}_{F_0}} {\mathcal O}_F$ and  to consider  the ring $\mathcal{O}_{F}[[u]][u^{-1}]$. The ring $\mathbb{A}_{F,\pi}$ will be the $\pi$-adic completion of this ring. The action of $\Gamma_{F}$ will be trivial on $\mathcal{O}_{F}$; on $u$ it is given by $\gamma u=[\chi_{F}](\gamma)(u)$, and the $r$-th iterate of the Frobenius is trivial on $\mathcal{O}_{F}$ and on $u$ is given by $\varphi^{r}(u)=[\pi](u)$. The residue field of $\mathbb{A}_{F,\pi}$ is $k_{F}((u))$ which is just $\mathbb{E}_{F,\pi}$, and the action of $\Gamma_{F}$ and $\varphi^{r}$ lifts the natural action of $\Gamma_{F}$ and $\varphi^{r}$ on $\mathbb{E}_{F,\pi}$.
~$\Box$\\ 
Let ${\mathbb A}_{F,\pi}^{ur}$ be the ring of integers of the unramified closure ${\mathbb B}^{ur}_{F,\pi}$ of the fraction field $\mathbb{B}_{F,\pi}$ of the ring $\mathbb{A}_{F,\pi}$.  We denote by ${\hat {\mathbb A}}_{F,\pi}^{ur}$ the completion of  $ \mathbb{A}^{ur}_{F,\pi}$.
\begin{lem}
The residue field of $\mathbb{A}_{F,\pi}^{ur}$ is a separable algebraic closure $\mathbb{E}_{\pi}$ of $\mathbb{E}_{F,\pi}$. Moreover, there is a natural isomorphism between the Galois group of ${\mathbb B}^{ur}_{F,\pi}$ over  ${\mathbb B}_{F,\pi}$ and $H_{F}$.

\end{lem}
{\bf Proof.}
See \cite{Ki-R}, lemma (1.4).
~$\Box$\\ 
We may now define the functors $V$ and $D_{\mathbb{A}_{F,\pi}}$. Let $A$ be a $p$-adic representation of $G_{F}$ over $\mathcal{O}_{F}$, 
then we define $D_{\mathbb{A}_{F,\pi}}(A):=(A\otimes_{\mathcal{O}_{F}}{\hat {\mathbb A}}_{F,\pi}^{ur})^{H_{F}}$.
 Let $B$ be an \'etale $(\varphi^{r},\Gamma_{F})$-module over $\mathbb{A}_{F,\pi}$, then we define $D_{\mathbb{A}_{F,\pi}}(B):=(B\otimes_{\mathbb{A}_{F,\pi}} {\hat {\mathbb A}}_{F,\pi}^{ur})^{\varphi^{r}=1}$.
\begin{thm} \label{AA}
The functors $V$ and $D_{\mathbb{A}_{F,\pi}}$ are quasi-inverse equivalences between the categories ${\bf{Rep}}_{G_{F}}(\mathcal{O}_{F})$ and ${\bf{\Phi}}^{et}_{\varphi^{r},\Gamma_{F}}({\mathbb{A}_{F,\pi}})$. 
\end{thm}
{\bf Proof.}
See \cite{Ki-R}, Theorem (1.6). 
$\Box$\\ 

\section{The ring $\mathbb{A}_{K,\pi}$ for $K$ a finite extension of $F$}

We now consider a finite extension $K$ of $F$:  unlike the situation considered by Kisin and Ren \cite{Ki-R} we place no  hypothesis 
on the extension .  Let $F_{\pi}$ be the maximal abelian ramified extension of $F$ attached to the uniformizer $\pi$ by Lubin-Tate theory. By \cite{W} 3.2, the field of norms $\mathbb{E}_{K,\pi}$ of the extension $K\cdot F_{\pi}$ of $K$ is a finite separable extension  of $\mathbb{E}_{F,\pi}$. Let $\mathbb{A}_{F,\pi}^{ur}$ be the strict Henselization of $\mathbb{A}_{F,\pi}$.  We remark again that the uniformizer of $\mathbb{A}_{F,\pi}$ is $\pi$.  There exists a unique subring of $\mathbb{A}_{F,\pi}^{ur}$, which is a discrete valuation ring  containing $\mathbb{A}_{F,\pi}$ and whose residue field is $\mathbb{E}_{K,\pi}$. Let us denote this ring by $\mathbb{A}_{F,\pi}(K)$. Our approach follows that given in \cite{Co} \S4,6.
\begin{pro}
The ring  $\mathbb{A}_{F,\pi}(K)$ is stable under $G_{K}$ and $\varphi$.
\end{pro}
{\bf Proof.}
On $\mathbb{A}_{F,\pi}^{ur}$, the action of $G_{F}$ lifts the action of $G_{F}$ on the fixed separable algebraic closure of $\mathbb{E}_{F,\pi}$, and the action of $\varphi$ lifts the action of the Frobenius. Since $G_{K}$ and $\varphi$ leave $\mathbb{E}_{K,\pi}$ stable, also $\mathbb{A}_{F,\pi}(K)$ will be stable under $G_{K}$ and $\varphi$.
~$\Box$\\ 

{\it Remark.}     Note that, if we put ourselves under the hypotheses of  \cite{Ki-R} 2.1 (which are in force throughout that article), 
our $\mathbb{A}_{F,\pi}(K)$ coincides with     $\mathcal{O}_{\mathcal{E}}$ of  \cite{Ki-R}.
\medskip

The ring   $\mathbb{A}_{F,\pi}(K) $  contains the ring $\mathcal{O}_{K_{0}}\otimes_{\mathcal{O}_{F_{0}}} \mathcal{O}_{F}$  (we recall that   $\mathcal{O}_{F_{0}}$ and $\mathcal{O}_{K_{0}}$ denote respectively the ring of Witt vectors on $k_{F}$ and $k_{K}$ ). 
The ring $\mathbb{A}_{F,\pi}(K)$ is unramified over the discrete valuation ring  $\mathcal{O}_{K_{0}}\otimes_{\mathcal{O}_{F_{0}}} \mathcal{O}_{F}$.  Moreover $\mathcal{O}_{K}$ is totally ramified over $\mathcal{O}_{K_{0}}\otimes_{\mathcal{O}_{F_{0}}} \mathcal{O}_{F}$. 
Hence the tensor product $\mathbb{A}_{F,\pi}(K) \otimes_{\mathcal{O}_{K_{0}}\otimes_{\mathcal{O}_{F_{0}}} \mathcal{O}_{F}}\mathcal{O}_{K}$ is again a discrete valuation ring.  Let us indicate it  by  $\mathbb{A}_{K ,\pi}$. It has uniformizer $\varpi$ and residue field $\mathbb{E}_{K,\pi}$. 
We summarize the properties of  $\mathbb{A}_{K ,\pi}$: it arises from the  field of norms of the   maximal totally ramified extension of $K$ related  to the  Lubin-Tate group on $F$  given by   $\pi$ and we have added to its Cohen lifting     the ramification of $K$.  Moreover, one may extend the  commuting action of $G_{K}$ and of $\varphi^{rs}$ on ${\mathbb A}_{F,\pi}(K)$ to (commuting) actions on ${\mathbb A}_{K,\pi}$  by imposing that the  $rs$-iterate of Frobenius  act trivially on $\mathcal{O}_{K}$.

We have just proved the following.
\begin{lem}
There exists a complete discrete valuation ring  $\mathbb{A}_{K,\pi}$, with uniformizer $\varpi$ and residue field $\mathbb{E}_{K,\pi}$, endowed with commuting actions of $G_{K}$ and $\varphi^{rs}$ lifting the natural actions of $G_{K}$ and $rs$-th iterate of the Frobenius on $\mathbb{E}_{K,\pi}$.
\end{lem}
{\bf Proof.}
The construction above. Observe that $K$ is  fixed  under the $rs$-th iterate of the Frobenius.
~$\Box$

We may then take the unramified closure of  $ Frac({\mathbb A}_{K,\pi})$ in  the extension  of $\widetilde {\mathbb B}$ by  $K$ (denote it by ${\widetilde {\mathbb B}}_K$). Note that  while  $K$ is totally ramified over the $K_0$,  ${\widetilde {\mathbb B}}$ is unramified over $K_0$.  We will indicate its valuation ring as
${\mathbb A}_{K,\pi} ^{ur}$ and by  ${\hat{\mathbb A}}_{K,\pi}
^{ur}$ its completion. We then have: 
\begin{lem}
The ring ${\hat{\mathbb A}}_{K,\pi}
^{ur}$ is a complete discrete valuation ring  with residue field a separable algebraic closure of $\mathbb{E}_{K,\pi}$; it is equipped with commuting actions of $G_{K}$ and $\varphi^{rs}$, lifting the actions on the separable algebraic closure $\mathbb{E}_{\pi}$ of $\mathbb{E}_{K,\pi}$. 
\end{lem}
{\bf Proof.}
The ring ${\mathbb A}_{K,\pi}
^{ur}$ may be viewed as
 $\mathbb{A}_{F,\pi}^{ur} \otimes_{\mathcal{O}_{K_{0}}\otimes_{\mathcal{O}_{F_{0}}} \mathcal{O}_{F}}\mathcal{O}_{K}$ 
 since $K$ is totally ramified over the $F \cdot K_0$ and $Frac({\mathbb A}_{F,\pi}(K))$ is unramified over $F \cdot K_0$. 
 The action of $G_{K}$ extends that on  $\mathbb{A}_{F,\pi}^{ur}$ by the trivial action on
  $\mathcal{O}_{K}$, and analogously for the action of $\varphi^{rs}$ in the completion:
   ${\hat{\mathbb A}}_{K,\pi}^{ur} ={\hat {\mathbb A}}_{F,\pi}^{ur} \otimes_{\mathcal{O}_{K_{0}}\otimes_{\mathcal{O}_{F_{0}}} \mathcal{O}_{F}}\mathcal{O}_{K}$.  ~$\Box$

\medskip

In particular $Gal ({\hat{\mathbb A}}_{K,\pi}
^{ur}/ {\mathbb A}_{K,\pi})= Gal ( {\mathbb E}_{\pi} /{\mathbb E}_{K, \pi})= H_{K,\pi}$.  From now on,  in order to lighten  the notation,  we will denote ${\hat{\mathbb A}}_{K,\pi}
^{ur} $ by ${\mathbb A}_{\pi}$.  Note that $G_K$ acts on ${\mathbb A}_{K,\pi}$ through $\Gamma_{K,\pi}=  \chi_{\pi}(G_K)$. We may now define the functors $V_{\pi}$ and $D_{\pi}$.
Let $A$ be $p$-adic representation of $G_{K}$ over $\mathcal{O}_{K}$, then we define $D_{\pi}(A):=(A\otimes_{\mathcal{O}_{K}} {\mathbb A}_{\pi})^{H_{K,\pi}}$.
 Let $B$ be an \'etale $(\varphi^{rs},\Gamma_{K})$-module over ${\mathbb A}_{K,\pi}$, 
 then we define $V_{\pi}(B)
:=(B\otimes_{\mathbb{A}_{K,\pi}}   {\mathbb A}_{\pi})^{\varphi^{rs}=1}$.
\begin{thm}\label{teo2}, 
The functors $V_{\pi}$ and $D_{\pi}$ are quasi-inverse equivalences between the categories ${\bf{Rep}}_{G_{K}}(\mathcal{O}_{K})$ and ${\bf{\Phi}}^{et}_{\varphi^{rs},\Gamma_{K}}(\mathbb{A}_{K,\pi})$. 
\end{thm}
{\bf Proof.}
To prove the statement we proceed as in \cite{Fo} 1.2.4, 1.2.6. 
It is enough to prove that the maps 
$$
 D_{\pi}(V) \otimes_{\mathbb{A}_{K,\pi}} {\mathbb A}_{\pi} 
\rightarrow 
 V   \otimes_{\mathcal{O}_{K}} {\mathbb A}_{\pi}
$$
and 
$$ V_{\pi}(B)  \otimes_{\mathcal{O}_{K}}  {\mathbb A}_{\pi}
\rightarrow  
B    \otimes_{\mathbb{A}_{K,\pi}} {\mathbb A}_{\pi}
 $$
are isomorphisms.
This follows in the limit from the analogous statement on torsion modules, \textit{i.e.}, modules killed by a power of $\varpi$. 
The statement on torsion modules follows by d\'evissage from the statement on modules killed by $\varpi$. On modules killed by $\varpi$, the statement was proved in the proof of Theorem \ref{teor1}.
~$\Box$\\

{\it Remark.}  As we mentioned the results presented in this paragraph are generalizations of the contructions given in \cite{Ki-R}.  
Throughout that article the two authors  made the  assumption  ( \cite {Ki-R}2.1) that,  if we denote  by $\{ F_n \} $ the totally ramified tower over $F$ associated to the Lubin-Tate group given by $\pi$ and by $K_{0,L}$ the field $K_0 \otimes _{L_0} L$ then there exists $n$ such that 
$K \subset F_n \cdot K_{0,L} $. This  condition will imply that the field of norms  ${\mathbb E}_{K,\pi}$ is unramified over ${\mathbb E}_{F, \pi}$ 
and the degree of the extension is  $[K_0:L_0]$.  Our general setting will allow ramification and 
the  residue degree  will be the degree  over $F$ of  the maximal unramified extension of $F$ in the 
maximal totally ramified extension of $K$ given by the union of the  $F_n \cdot K$.

\section{Comparison between the two different equivalences}

In the preceding two sections we have seen two different equivalences involving the category ${\bf {Rep}}_{G_{K}}(\mathcal{O}_{K})$ and two categories of $(\varphi,\Gamma)$-modules. 
The first equivalence is obtained via the method explained in Section 6, hence we  start with a Lubin-Tate group on $F$, associated to a uniformizer $\pi$  and  then we can consider the totally ramified tower attached to it over $K$.  
But we could have started directly from $K$, chosen a uniformizer $\varpi$ in it, and then  built another totally ramified tower: then 
we would be in the case of Section 5 (but with $F$ replaced by $K$ in that section).  
 In the notation of Section 5, then we would then have an equivalence between  $\bf{Rep}_{G_{K}}(\mathcal{O}_{K})$ and $(\varphi^{rs}, \Gamma_{K,\varpi})$-modules over 
  ${\mathbb A}_{K,\varpi}$.   Theorem \ref{AA}  now would read as follows.
\begin{thm}

The two functors $$D_{\varpi}:{\bf{Rep}}_{G_{K}}(\mathcal{O}_{K})\rightarrow {\bf{\Phi}}^{et}_{\varphi^{rs},\Gamma_{K,\varpi}}(\mathbb{A}_{K,\varpi})$$
 and 
$$V_{\varpi}:{\bf{\Phi}}^{et}_{\varphi^{rs},\Gamma_{K,\varpi}}(\mathbb{A}_{K,\varpi}) \rightarrow {\bf{Rep}}_{G_{K}}(\mathcal{O}_{K}),$$
defined 
by $D_{\varpi}(A)=(A\otimes_{\mathcal{O}_{K}}{\hat{\mathbb A}}_{K,\varpi}^{ur} )^{H_{K,\varpi}}$ 
for
  $A\in {\bf{Rep}}_{G_{K}}(\mathcal{O}_{K})$, 
  and $V_{\varpi}(B)=(  B \otimes_{\mathbb{A}_{K,\varpi}}{\hat{\mathbb A}}_{K,\varpi}^{ur})^{\varphi^{rs}=1}$, 
  for $B\in {\bf{\Phi}}^{et}_{\varphi^{rs},\Gamma_{K,\varpi}}({\mathbb A}_{K,\varpi} )$, induce quasi-inverse equivalences of categories.
\end{thm}
 Again in order to  lighten the notation   we will indicate  ${\hat{\mathbb A}}_{K,\varpi}^{ur}$ by 
 ${\mathbb A}_{\varpi}$.  So we have that the category ${\bf{\Phi}}^{et}_{\varphi^{rs},\Gamma_{K,\varpi}}(\mathbb{A}_{K,\varpi})$ is equivalent to ${\bf{Rep}}_{G_{K}}(\mathcal{O}_{K})$. But, by Theorem \ref{teo2},  the category ${\bf{\Phi}}^{et}_{\varphi^{rs},\Gamma_{K,\pi}}(\mathbb{A}_{K,\pi})$ is equivalent to ${\bf{Rep}}_{G_{K}}(\mathcal{O}_{K})$. By transitivity, we have that the category ${\bf{\Phi}}^{et}_{\varphi^{rs},\Gamma_{K,\varpi}}(\mathbb{A}_{K,\varpi})$ is equivalent to  ${\bf{\Phi}}^{et}_{\varphi^{rs},\Gamma_{K,\pi}}(\mathbb{A}_{K,\pi})$.  Recall again that the first category uses the Lubin-Tate group defined by the uniformizer $\varpi$ of $K$, while the second uses the extension to $K$  of the Lubin-Tate group of $F$ associated  to $\pi$.

We wish to describe the two functors which give  the equivalence between these two categories  without appealing to the category ${\bf{Rep}}_{G_{K}}(\mathcal{O}_{K})$.


 First  we define a system of rings $\mathbb{A}_{\pi,\varpi,i}$, together with 
 a $G_{K}$-action on each of these rings, two partial Frobenius morphisms $\varphi_{\pi}:\mathbb{A}_{\pi,\varpi,i}\rightarrow\mathbb{A}_{\pi,\varpi,i-1}$ and 
 $\varphi_{\varpi}:\mathbb{A}_{\pi,\varpi,i}\rightarrow\mathbb{A}_{\pi,\varpi,i+1}$. 
 These rings will be complete discretely valued rings with uniformizer $\varpi$ and the residue field of $\mathbb{A}_{\pi, \varpi,i}$ being
  $\mathbb{E}_{\pi, \varpi,i}$; the $G_{K}$-action and the partial Frobenius morphisms are compatible with reduction mod $\varpi$. 
  The system of fields $\mathbb{B}_{\pi,\varpi,i}$ is obtained from the system of rings $\mathbb{A}_{\pi,\varpi,i}$ by inverting $\varpi$; of course it inherits the $G_{K}$-action and the partial Frobenius morphisms.

Let us now define the rings $\mathbb{A}_{\pi,\varpi,i}$.
Consider the rings ${\mathbb A}_{\pi}= {\hat {\mathbb{A}}}_{K,\pi}^{ur}$ 
 and 
 ${\mathbb A}_{\varpi}= {\hat {\mathbb{A}}}_{K,\varpi}^{ur}$  of the previous sections.  Consider the tensor product 
$ {\mathbb A}_{\pi} \otimes_{W(\overline{k}_{K}),\sigma^{i}} {\mathbb A}_{\varpi}$, 
where $W(\overline{k}_{K})$ is the ring of Witt vectors with coefficients in the field $\overline{k}_{K}$, and the tensor product is 
twisted by the $rsi$-th power of the Frobenius automorphism of $W(\overline{k}_{K})$. 
In this tensor product, the ideal $\varpi(  {\mathbb A}_{\pi} \otimes_{W(\overline{k}_{K}),\sigma^{i}} {\mathbb A}_{\varpi})$ is prime, since the quotient is $\mathbb{E}_{\pi}\otimes_{\overline{k}_{K},\sigma^{i}}\mathbb{E}_{\varpi}$, which is a domain 
by Lemma \ref{lem2}.    Let $T$ be the localization of $  {\mathbb A}_{\pi} \otimes_{W(\overline{k}_{K}),\sigma^{i}} {\mathbb A}_{\varpi}$  with respect to the prime ideal
 $\varpi( {\mathbb A}_{\pi} \otimes_{W(\overline{k}_{K}),\sigma^{i}} {\mathbb A}_{\varpi})$. 
 It is then a local ring with residue field $\mathbb{E}_{\pi,\varpi,i}$ and its maximal ideal is generated by
  $\varpi$. The ring $\mathbb{A}_{\pi,\varpi,i}$ is defined as the $\varpi$-adic completion of the ring $T$; and $\mathbb{B}_{\pi, \varpi,i}$ is the quotient field $\mathbb{A}_{\pi,\varpi,i}[1/\varpi]$ of $\mathbb{A}_{\pi,\varpi,i}$.

The $G_{K}$-action and the partial Frobenius maps of ${\mathbb A}_{\pi} \otimes_{W(\overline{k}_{K}),\sigma^{i}} {\mathbb A}_{\varpi}$ induce a $G_{K}$-action and partial Frobenius morphisms on the system $\mathbb{A}_{\pi,\varpi,i}$, and hence also on the system $\mathbb{B}_{\pi, \varpi,i}$. Clearly, on 
reducing $\mathbb{A}_{\pi,\varpi,i}$   mod $\varpi$  one obtains the system $\mathbb{E}_{\pi,\varpi,i}$.

The following proposition is the lifting of Proposition \ref{pro1} to characteristic 0.
\begin{pro}\label{pro2}
The projective limit of the system $\mathbb{A}_{\pi,\varpi,\bullet}$ along $\varphi_{\pi}$ is ${\mathbb{A}}_{ \varpi}$, and the morphisms $\varphi_{\varpi}$ are transformed into the Frobenius morphism $\varphi_{\varpi}$ of ${\mathbb{A}}_{\varpi}$. The symmetrical statement holds also and the same holds true for the system  $\mathbb{B}_{\pi,\varpi,\bullet}$.
\end{pro}
{\bf Proof.}
An element $(x_{i})$ is in the projective limit of the system $\mathbb{B}_{\pi,\varpi,i}$ for the morphism $\varphi_{\pi}$ if and only if, there exists a $t\in\mathbb{N}$, such that $(\varpi^{t}x_{i})$ is an element of the projective limit of the system $\mathbb{A}_{\pi,\varpi,i}$ for the morphism $\varphi_{\pi}$. So one is reduced to the computation of the  projective limit of the system $\mathbb{A}_{\pi,\varpi,i}$ for the morphism $\varphi_{\pi}$. Of course, any element of the type $(x_{i})$, with $x_{i}=x\in\mathbb{A}_{\varpi}$ is in this projective limit; so one must just show that any $(x_{i})$ in this projective limit is of this type.

Consider the reduction mod $\varpi$ of the element $(x_{i})$, it is an element of the projective limit of the system $\mathbb{E}_{\pi,\varpi,i}$ along the morphism $\varphi_{\pi}$; so 
by Proposition \ref{pro1}, there exists an element 
$x^{(0)}\in {\mathbb{A}}_{\varpi}$ such that for every $i$, $x_{i}\cong x^{(0)}$ mod $\varpi$.
 It follows that $(x_{i}-x^{(0)})$ is an element of the projective limit divisible by $\varpi$; so also $(x_{i}-x^{(0)})/\varpi$
  is an element of the projective limit. One may apply the same argument to get an element $x^{(1)}\in\mathbb{A}_{\varpi}$ such that, for every $i$, $(x_{i}-x^{(0)})/\varpi$ is congruent to $x^{(1)}$ mod $\varpi$. So $x^{(0)}+\varpi x^{(1)}$ is 
  congruent to $x_{i}$ mod $\varpi ^{2}$. By iterating this process,  for every index $i$ one obtains a convergent sequence of elements 
  of ${\mathbb{A}}_{\varpi}$, converging to $x_{i}$ . 
  So this projective limit can be identified with $\mathbb{A}_{\varpi}$; 
  and in this way, the morphism $\varphi_{\varpi}$ is identified with the usual Frobenius morphism  $\varphi^{rs}$ of $\mathbb{A}_{\varpi}$.

The symmetrical statement is dealt with symmetrically.
$\Box$\\ 




We can now state our main result:

\begin{thm}
There is an equivalence of categories between the categories ${\bf{\Phi}}^{et}_{\varphi^{rs},\Gamma_{K,\pi}}({\mathbb{A}_{K,\pi}})$ and ${\bf{\Phi}}^{et}_{\varphi^{rs},\Gamma_{K,\varpi}}({\mathbb{A}_{K,\varpi}})$. This equivalence is realized by the following functors

$$D_{1}\mapsto (\lim_{\leftarrow,\varphi_{\pi}}D_{1}\otimes _{\mathbb{A}_{\varpi,\pi}}\mathbb{A}_{\pi,\varpi,\bullet})^{H_{K,\varpi}};$$

and

$$D_{2}\mapsto (\lim_{\leftarrow,\varphi_{\varpi}}D_{2}\otimes _{\mathbb{A}_{K,\varpi}}\mathbb{A}_{\pi,\varpi,\bullet})^{H_{K,\pi}}.$$

\end{thm}
{\bf Proof.}
Consider the \'etale $(\varphi^{rs},\Gamma_{K,\pi})$-module $D_{1}$ over the ring $\mathbb{A}_{K,\pi}$, coming from the representation $V$ of $G_{K}$ over $\mathcal{O}_{K}$. 
Then we have,
 $D_{1}\otimes _{\mathbb{A}_{K,\pi}}\mathbb{A}_{\pi,\varpi,\bullet}=V\otimes _{\mathcal{O}_{K}}\mathbb{A}_{\pi,\varpi,\bullet}$. 
 Hence, 
 $(\lim_{\leftarrow,\varphi_{\pi}}D_{1}\otimes _{\mathbb{A}_{K,\pi}} \mathbb{A}_{\pi,\varpi,\bullet}=
 V\otimes_{\mathcal{O}_{K}}  \lim_{\leftarrow, \varphi_{\pi}}\mathbb{A}_{\pi,\varpi,\bullet}$. 
 By Proposition \ref{pro2}, the projective limit is $ {\mathbb{A}}_{\varpi}$. Thus $((\lim_{\leftarrow,\varphi_{\pi}}D_{1}\otimes _{\mathbb{A}_{K,\pi}} \mathbb{A}_{\pi,\varpi,\bullet})^{H_{K,\varpi}}=(V\otimes_{\mathcal{O}_{K}}  {\mathbb{A}_{\varpi}})^{H_{K,\varpi}}$. This shows that it is the equivalence of categories obtained by composing the two equivalences of categories with the category of representations of $G_{K}$ over the field $\mathcal{O}_{K}$.

The second statement is dealt with analogously. 
$\Box$\\

\bigskip

Bruno Chiarellotto, Dip. Matematica, Univ. Padova, Via Trieste 63, 35121 Padova (Italy), email: chiarbru@math.unipd.it
	
Francesco Esposito, Dip. Matematica, Univ. Padova, Via Trieste 63, 35121 Padova (Italy), email: esposito@math.unipd.it

\begin{thebibliography}{10}

\bibitem{Co}
P. Colmez, "Espace Vectoriels de dimension finie et representation de deRham", Asterisque 319 (2008), pp. 117-186.


\bibitem{Fo}
J-M. Fontaine, "Representations $p$-adiques des corps locaux I", in The Grothendieck festschrift, vol. II, 249-309, Progr. Math. 87, Birkhauser Boston, 1990.

\bibitem{FW}
J-M. Fontaine, J-P. Wintenberger, "Le corps des normes de certaines extensions algebriques de corps locaux", C.R.A.S. 288, 367-370 (1979).

\bibitem{Fou}
M.L. Fourquaux, " Logarithme de Perrin-Riou pour des extensions associ\'ees \`a un groupe de Lubin-Tate", These, Paris VI, 2005. Avalaible at http://www.normalesup.org/~fourquau/pro/publications/ 

\bibitem{Ha} 
M.Hazewinkel  "Formal groups and applications". Academic Press 1978.

\bibitem{Ki-R}
M. Kisin, W. Ren, "Galois representations and Lubin-Tate groups", Documenta Math. 14 (2009), pp.441-461.


\bibitem{LT}
J. Lubin, J. Tate, "Formal complex multiplication in local fields", Ann. of Math. 81, 380-387 (1965).

\bibitem{Se}
J.P. Serre, "Local classfield theory", in Cassels, Frolich, Algebraic number theory, Academic press, 1967.


\bibitem{W}
J.P. Wintenberger, "Le corps des normes de certaines extensions infinies de corps locaux; applications", Ann. scient. E.N.S.16, 59-89 (1983).



\end{thebibliography}
\end{document}